\documentclass[10pt]{article}
\usepackage{amsthm, amssymb,amsthm, latexsym,epsfig,color, 
graphicx, cite, setspace, enumitem, hyperref}
\usepackage[misc,geometry]{ifsym} 
\usepackage[normalem]{ulem}
\usepackage{tikz}
\usetikzlibrary{arrows}
\usepackage{tkz-euclide}
\usetkzobj{all}
\usepackage{extarrows, url}
\usepackage{authblk}

\newcommand{\R}{{\mathbb R}}

\newcommand{\G}{{\mathcal G}}

\newtheorem*{theorem*}{Theorem} 

\title
{On the equivalence of Playfair's axiom \\ to the parallel postulate 
}

\author[ ]{Elizabeth T. Brown\thanks{Department of Mathematics \& Statistics, James Madison University, Harrisonburg VA 22807. ({\tt brownet@jmu.edu}) } }
\author[ ]{Emily Castner
\thanks{Department of Mathematics \& Statistics, Mount Holyoke College, South Hadley, MA 01075.  
({\tt castn22e@mtholyoke.edu}) }
}
\author[ ]{Stephen Davis
\thanks{Department of Mathematics, University of Wisconsin, Madison, WI 53706. 
({\tt smdavis7@wisc.edu}) }
}
\author[ ]{Edwin O'Shea (\Letter)
\thanks{Department of Mathematics \& Statistics, James Madison University, Harrisonburg VA 22807. 
({\tt osheaem@jmu.edu}) }
}
\author[ ]{Edouard Seryozhenkov\thanks{Department of Mathematics, Willamette University, Salem, OR 97301. 
({\tt edseryo@gmail.com}) } 
}
\author[ ]{AJ Vargas\thanks{Department of Mathematics, Bryn Mawr College, Bryn Mawr, PA 19010. 
({\tt avargasjr@brynmawr.edu}) } \thanks{The authors received support from NSF grant DMS-1560151.}
}

\affil[ ]{}

\begin{document}

\maketitle

\begin{abstract}
We show that the classical equivalence of Euclid's parallel postulate and Playfair's axiom 
collapses in the absence of triangle congruence. 
In particular, we construct a non-SAS geometry that models the Playfair axiom 
but not the parallel postulate. 
\end{abstract}


\section*{Introduction}

Euclid's account of the Side-Angle-Side condition for triangle congruence, \cite[Proposition 4 of Book I]{Euc02}, 
relies on an appeal to superposition, the application of one triangle to another by  ``dragging-and-dropping.'' 
Yet he is at pains to avoid this method in Proposition 2  
-- and in many propositions that follow -- 
which would be easy to prove via the superposition of one line segment or angle to another. 
Indeed his use of superposition is conservative, used only in Propositions 4 and 8 with latter propositions 
like Proposition 23 rhyming with Proposition 2 in their restraint.

It is now understood that the presence of the SAS axiom in a geometry is equivalent to that geometry having 
rigid motions (cf. \cite[\S 17]{Har13}), so the preservation of distance and angles 
under motion is equipotent to the ability to declare two triangles congruent.
This is but one instance of Klein's ``Erlanger Programm'' which proposed that a geometry be completely classified by 
the invariant theory of its group of transformations \cite[Part III, \S 1.1]{Kle39}. 
This was a paradigm shift that had as pervasive an influence on mathematics as the discovery of non-Euclidean geometry. 
In the Euclidean context, Klein's approach is inescapably tied to 
triangle congruence, motivating our study of geometries borne of 
Euclidean axioms but without SAS.  

Indeed the failure of SAS has profound consequences for the parallel postulate. Classically, Euclid's parallel 
postulate has many equivalents, most notably Playfair's axiom which was known 
since antiquity \cite{Pro92} (cf. \cite[p.220]{Hea56}). 
The main result of this paper is as follows:  

\begin{theorem*} \label{the:main}
In the absence of SAS, Euclid's parallel postulate is not equivalent to Playfair's axiom. 
In particular, there exists a non-SAS geometry $\G$ that models the Playfair axiom 
but not the parallel postulate. 
\end{theorem*}


So far as we know, this is the first demonstration of this result. 
Beeson, in his 2016 monograph ``Constructive Geometry and the Parallel Postulate," 
\cite{Bee16} shows that the Playfair axiom is 
strictly weaker than the Parallel Postulate in the setting of constructive logic. His result differs 
from ours in that the rules of inference in constructive logic are restricted from those of classical 
first order logic, which is to say, from the rules of inference assumed by most mathematicians. 
For example, in constructive logic if one shows that the statement ``A does not exist" leads to a 
contradiction, one cannot conclude (without additional argument) that therefore A {\em{does}} exist. 
Moreover, Beeson's models are ``given without discussing the 
exact choice of the other axioms of geometry'' and thus are not concerned 
with non-SAS~\cite[Introduction]{Bee16}.

\section*{Axiomatic Context}

In keeping with modern sensibilities, we will use Hilbert's framework 
for Euclidean geometry vis-\`{a}-vis 
{\it Foundations of Geometry} \cite[Chapter I]{Hil71}. 
His axioms are grouped according to 
{\it incidence in the plane} (Axioms I.1-3), 
{\it order of points} or {\it betweeness} (Axioms II.1-4), 
{\it congruence for segments, angles, and triangles} (Axioms III.1-5), 
and the {\it axiom of parallels} (Axiom IV). 
For reference, we list the central axioms of SAS and parallels. 
Hilbert's version of the SAS axiom (III. 5) is:    

\vspace{.1in}

\parbox{.5in}
{{\bf (SAS)}} 
\hspace{.1in}
\parbox{5in}
{If two triangles share two sides and the angle in between those two sides then 
in each triangle there are second angles that are also equal.
}

\vspace{.1in}

\noindent This says nothing of the third side and third angle but Hilbert 
shows \cite[\S 6, Theorem 12]{Hil71} that (classical) Euclidean 
SAS  holds assuming axiom groups I, II, and III. 
Hilbert's axiom of parallels, Axiom IV \cite[\S 4]{Hil71}, curiously called 
``Euclid's Axiom'' by Hilbert, states: 

\vspace{.1in}

\parbox{.5in}
{{\bf (hPF)}} 
\hspace{.1in}
\parbox{5in}
{Let $a$ be any line and $A$ a point not on it in a common plane. 
Then there is {\it at most one line} in the plane, determined by $a$ and $A$, 
that passes through $A$ and does not intersect $a$.
}

\vspace{.1in}

\noindent Hilbert's version is slightly weaker than the classical Playfair axiom 
{\bf (cPF)}, which insists that there is \textit{exactly one line} 
rather than merely \textit{at most one line}.  
Hilbert's version allows for, say, the geometry of geodesic lines on the sphere. 
Euclid's original parallel postulate \cite[Book I, Postulates]{Euc02} 
asserts: 
 
\vspace{.1in}

\parbox{.5in}
{{\bf (PP)}} 
\hspace{.1in}
\parbox{5in}
{That, if  a straight line falling on two straight lines make the 
interior angles on the same side less than two right angles, the two 
straight lines, if produced indefinitely, meet 
on that side on which are the angles less than the two right angles. 
}

\vspace{.1in}

The reader will note that Euclid's statement involves both lines and angles, 
his two basic objects of plane geometry, while Playfair refers to lines only. 
Most authors of the late 19th century chose Playfair over the parallel postulate, 
a tradition that persists in modern treatments of the foundations of Euclidean 
geometry like Hartshorne \cite{Har13}. 
Dodgson \cite{Dod85} is the most notable holdout, exalting the positive 
constructive formulation 
of Euclid's original postulate over the existence assertion of Playfair 
(cf. \cite[p.313-4]{Hea56}). It is not obvious from the literature why this 
axiomatic drift to Playfair's version occurred. In Hilbert's case 
it seems reasonable that the success of his non-constructive approach 
to invariant theory (see \cite[Pages {\it xix} and 40]{Olv99}) 
had a significant influence on how he approached axiomatics in geometry. 
Or he may simply have been scrupulous in applying the aesthetic of parsimony. 

\section*{Proof of the Main Theorem}

We will construct a geometry $\G$ that respects Hilbert's
axioms of incidence (I), order (II), and congruence (III) except 
for the SAS axiom (III.5). The geometry models cPF, and therefore hPF, 
but not PP. Points and lines of $\G$ will be those of the usual 
Euclidean plane $\R^2$. Angle congruence is manipulated 
from the usual measure while still maintaining the literal 
integrity of Hilbert's axioms.

Hilbert defines angles in terms of rays.  
Given a line $\bar{h}$ and a point $O$ on $\bar{h}$, 
there are two {\it rays} (or ``half lines''), $h$ and $h^\prime$ (well defined by the 
group of order axioms, Group II) on either side of $O$. 

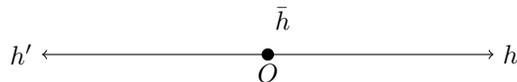
\begin{center}
\begin{tikzpicture}
    \draw[<-] (-3,0) -- (0,0);
    \draw[->] (0,0) -- (3,0);
	\node [above] at (0.2,0.2) {$\bar{h}$};
	\node [left] at (-3,0) {$h^\prime$};
	\node [right] at (3,0) {$h$};
	\draw[fill] (0,0) circle [radius=0.075];
	\node [below] at (0,0) {$O$};
\end{tikzpicture}
\end{center}

\noindent
An {\em{angle}}
\footnote{By contrast, Euclid's Definition 8 is: 
	A {\em plane angle} is the inclination to one another of two lines in a plane 
	which meet one another and do not lie in a straight line. \cite{Euc02} 
	} 
is then defined as any two distinct rays 
$h, k$ with common vertex $O$ 
and lying on distinct lines $\bar{h}, \bar{k}$ respectively. 
The pair of rays $h,k$ is called an angle and is denoted by 
$\angle (h,k)$ or by $\angle (k,h)$ \cite[p.11]{Hil71}.

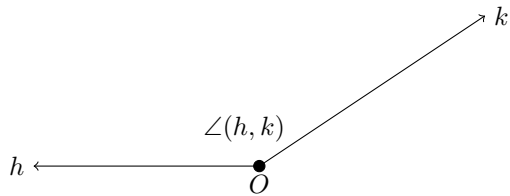
\begin{center}
\begin{tikzpicture}
    \draw[<-] (-3,0) -- (0,0);
    \draw[->] (0,0) -- (3,2);
	\node [above] at (-0.2, 0.2) {$\angle(h,k)$};
	\node [left] at (-3,0) {$h$};
	\node [right] at (3,2) {$k$};
	\draw[fill] (0,0) circle [radius=0.075];
	\node [below] at (0,0) {$O$};
\end{tikzpicture}
\end{center}

\noindent Hilbert excludes angles comprised of pairs of rays which together 
make a line. 
His definition of the interior of an angle is complicated by a desire to refer only to sides of rays and lines; he 
then uses axioms of incidence and order to prove that his definition selects the part of the plane that includes 
any line segment connecting the rays~\cite[p.11]{Hil71}. Since order of rays is irrelevant and interior 
of angle is unambiguous, every angle would have, in standard terms, radian measure exceeding $0$ but less then $\pi$.
Two angles having a vertex and one side in common and whose 
separate sides form a line are called {\em supplementary angles}; 
angles which are congruent to one of their supplementary 
angles are called {\em right angles}~\cite[p.13]{Hil71} 
There is no further restriction on congruence other than 
that every angle be congruent to itself~\cite[p.12]{Hil71}.

The relation of angle congruence in $\G$ will be interpreted as follows. 
For angle $\angle (h, k),$ let $\mu (\angle (h, k))$ be the standard radian measure without regard to order, 
so that $\mu (\angle (h, k)) = \mu (\angle (k, h))$. 
The range of $\mu$ on angles in this context is  $(0, \pi)$. 
For each point $P$ in $\R^2$, choose a bijection
$$
{\hat f}_P: (0, {\tfrac \pi 2}) \to (0, {\tfrac \pi 2}).
$$	
	\noindent
For each point $P$, we use the function ${\hat f}_P$ 
to label angles with vertex 
at $P$ in a way specific to $P$. For $\angle (h,k)$  with 
vertex $P$, define 
	\[f_P (\angle (h,k)) = \begin{cases}
	{\hat f}_P \left( \mu  (\angle (h,k))\right)  & {\rm{if}}\  \mu (\angle (h, k)) < {\tfrac \pi 2}\\
	{\tfrac \pi 2 } &  {\rm{if}}\  \mu (\angle (h, k)) = {\tfrac \pi 2}\\
	\pi - {\hat f}_P \left(\pi - \mu  (\angle (h, k))   \right) & {\rm{if}}\  \mu (\angle (h, k)) > {\tfrac \pi 2} .
	\end{cases}
	\]
	%

\noindent 
We interpret the congruence relation on angles using vertex labels. 
Given $\angle (h,k)$ with vertex $P$ and $\angle (h',k')$ with vertex $P'$,
\[ \angle (h,k)  \sim \angle (h',k') 
\, \, \,  
{\hbox{iff}} 
\, \, \,  
f_P\left( \angle (h,k) \right)   
=  f_{P'}\left( \angle (h',k') \right). \]
\noindent
Note that all right angles are in the same equivalence class, 
since $f_P ({\frac \pi 2}) = {\frac \pi 2}$ for all $P$.  
``Adding angles'' is accomplished by adding 
equivalence classes in the usual way. 
Supplementarity is also preserved: 
if two angles $\angle (h,k)$ and $\angle (k,l)$ 
with common vertex $P$ together make a line, then  
$ f_P\left( \angle (h,k) \right) + 
f_P\left( \angle (k,l) \right)
= \pi$. 

This relation respects 
Hilbert's angle congruence axiom, III.4: Each 
$\angle (h,k)$ is equivalent to itself. Since $f_P$ 
and $f_{P^\prime}$ are bijections, given a fixed $h,k$ with 
common vertex $P$ and a fixed $h'$ with vertex $P'$ and a given 
side of the line that contains $h'$, there is a 
unique $k'$ with vertex at $P'$ on the given side of 
the line that contains $h'$ such that 
$f_P\left( \angle (h,k) \right)   
=  f_{P'}\left( \angle (h',k') \right)$. 
Since Hilbert's axioms before III.4 do not concern 
angles -- only points and lines -- the geometry 
$\G$ satisfies Hilbert's axiom groups I, II, and III 
with the exception of III.5, the SAS axiom. SAS fails in $\G$ 
since rigid motions of a triangle do not preserve the angles 
of that triangle. 

To complete the construction fix values of the $f_P$ functions at 
the origin and at the point $(1,1)$, 
with lines $\ell$ given by $y = 0$, $\ell'$ given by $y = 1$, 
transversal $t$ given by $y = x$, 
and angles $\theta, \theta', \theta^{\prime\prime}$ as below.
 
 
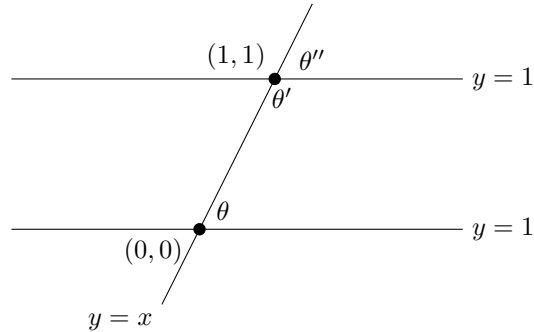
\begin{center}
\begin{tikzpicture}
  
    \draw (5,2) -- (11,2);
    \draw (5,0) -- (11,0);
    \draw (7,-1) -- (9,3);
	\node [right] at (11,2) {$y=1$};
	\node [right] at (11,0) {$y=1$};
	\node [below left] at (7,-1) {$y=x$};
	\node [above right] at (7.6,0) {$\theta$};
	\node [below left] at (7.4,0) {$(0,0)$};
	\draw[fill] (7.5,0) circle [radius=0.075];
	\node [below] at (8.6,2) {$\theta^\prime$};
	\node [above] at (9,2) {$\theta^{\prime\prime}$};
	\node [above left] at (8.5,2) {$(1,1)$};
	\draw[fill] (8.5,2) circle [radius=0.075];
\end{tikzpicture}
\end{center}
 
 \noindent
Fix bijections $f_{(0,0)}$ and $f_{(1,1)}$ so that 
$\displaystyle { f}_{(0,0)} \left({\tfrac \pi 4}\right) = {\tfrac \pi 4} $ 
and
$\displaystyle { f}_{(1,1)} \left({\tfrac \pi 4}\right) = {\tfrac {7 \pi} {16}} $. 
Then 
$\displaystyle \theta \in [{\tfrac \pi 4} ]$, 
and 
$\displaystyle \theta^{\prime\prime} \in [{\tfrac {7 \pi} {16}} ]$. 
Then $\theta^\prime$, $\theta^{\prime\prime}$'s supplement, satisfies 
$\displaystyle \theta^\prime \in [{\tfrac {9 \pi }{16}}]$. 
So 
$\displaystyle [\theta ] + [\theta^\prime] = [{\tfrac {13 \pi} {16}}]$, 
which is less than two right angles. 
The geometry $\G$'s angle congruency $\sim$ has not changed lines or intersections, 
so Playfair's axiom is still satisfied. But the parallel postulate fails, since 
the congruency class  $[\theta] + [\theta^\prime]$ is less than two right angles, 
but lines $y = 0$ and $y = 1$ do not intersect. This completes the proof. 

Although the assignments $f_{(0,0)}$ and $f_{(1,1)}$ were specific, the 
$f_P$ functions are quite flexible. To create 
a non-SAS geometry with this labeling scheme requires only one point $P$ at 
which angle congruence classes are assigned in an unorthodox but axiomatically legitimate way. 
It is also possible to insist that $\G$
preserve order on angles; to each point $P$ in the plane at which nonstandard labeling is desired, 
assign a distinct $r \in (1, \infty)$. Then let 
	\[{\hat f}_P (x) = x \left({\tfrac   {2x} \pi  }\right)^r  \qquad {\hbox { for } } \  x \in (0, {\tfrac \pi 2})\]
\noindent This assignment for every point $P$ in $\G$ would preserve Common Notion 5 of Euclid, 
that ``the whole is greater than the part'' and the manner in which it is invoked in 
the angle ordering of Propositions 16 and 20 of Book I. 

\section*{Closing Remarks}

The question of the converse is already settled. Hilbert (and others) showed that in the presence of all his axioms, 
Playfair and the parallel postulate are equivalent. So, a model of PF and $\lnot$PP must fail at least one 
of the other axioms. Without further ado, then, we can conclude that, if Hilbert's other axioms are assumed, 
$$
(\textup{PF} \land \lnot \textup{PP}) \rightarrow \lnot \textup{SAS}
$$ 
by Hilbert's own arguments. Our contribution lies 
in showing that the statement is not vacuously true since the antecedent is satisfiable.

\vspace{.3cm}

Playfair is the most-used but by no means only equivalent of the parallel 
postulate (cf. \cite[p.220]{Hea56}).  
The behavior of these equivalencies in absolute geometry (see Pambuccian ~\cite{Pam09})  
and other weakened forms of Euclidean geometry is of ongoing interest. 
For example, it would be natural to suppose that in our geometry $\G$ 
those equivalences which 
involve angles ought not to hold. This is false. 
Legendre's equivalence, ``there exists one triangle whose angle sum is two right angles,'' 
is based on angles only but $\G$ has enough flexibility to model Legendre's condition: 
Take any three points $Q,R, S$ and let $f_Q = f_R = f_S = \textup{id}$, 
the identity map. 
Then Legendre's condition holds on the triangle $\Delta QRS$ in $\G$. 

We close by asking for statements weaker than, or simply other than, SAS  
that are strong enough to recover all of the classical equivalencies of the parallel postulate.
For example, an intuitive property of angles that holds in the standard model of Euclidean geometry yet is
stronger than Common Notion 5, is that a whole angle ought to be the sum of its parts. 
Namely, given rays $h,k,r$ with common vertex and with $r$ ``inside'' the angle $\angle(h,k)$ 
then $\angle(h,k) = \angle(h,r) + \angle(r,k)$. 
Does a non-SAS geometry with this provision necessarily satisfy all the traditional equivalancies 
of PP? If so, does this provision form an intermediate axiom between non-SAS and SAS in the presence 
of Hilbert's other axioms?

\end{document}